\documentclass[11pt]{article}
\usepackage{amsmath,amssymb}


\begin{document}
\thispagestyle{empty} \setcounter{page}{1}

\noindent
%{\footnotesize {\bf To appear in\\[-1.75mm] {\em Dynamics of Continuous, Discrete and Impulsive Systems}}\\[-1.00mm]

%http:monotone.uwaterloo.ca/$\sim$journal} $~$ \\ [.3in]

%%USE THE STUFF BELOW AS A GUIDE TO SET UP THE START OF PAPER%%
%{\large{\bf { \begin{center}  OPTIMAL FEEDBACK CONTROL  FOR   DIFFERENTIAL INCLUSIONS ON  BANACH SPACES
%\end{center}}}
%\begin{center}  N.U.Ahmed \\ University of Ottawa \end{center}\vskip30pt
%
%
%%\begin{center}N.U.Ahmed,\\ University of Ottawa \end{center}\vskip30pt \begin{center}AMS %Joint Mathematics Meetings,
%%\\ SAN FRANCISCO,Jan 13-16, 2010\end{center}
%
%

{\large{\bf { \begin{center}   STOCHASTIC MINIMUM  PRINCIPLE FOR PARTIALLY OBSERVED SYSTEMS SUBJECT TO CONTINUOUS AND JUMP DIFFUSION PROCESSES  AND  DRIVEN BY  RELAXED CONTROLS
\end{center}}}

\vskip20pt

\begin{center}  N.U. Ahmed$^\dagger$ and C.D. Charalambous$^\ddagger$ \end{center}

\begin{center}
  University of Ottawa$^\dagger$ , Ottawa, Canada    \\ University of Cyprus$^\ddagger$, Nicosia, Cyprus. \end{center}

  \vskip30pt

%\begin{center}N.U.Ahmed,\\ University of Ottawa \end{center}\vskip30pt \begin{center}AMS %Joint Mathematics Meetings,
%\\ SAN FRANCISCO,Jan 13-16, 2010\end{center}

 \large \begin{center} {\bf  ABSTRACT }\end{center} In this paper we   consider non convex control problems of    stochastic differential equations driven by relaxed controls. We present existence of optimal controls and then develop  necessary conditions of optimality. We cover both continuous diffusion and Jump processes.


  \section{Introduction}
\label{introduction}

\par The basic idea of the deterministic minimum principle introduced by Pontryagin and his colleagues in the 1950's is to derive a set of necessary and sufficient  conditions that must be satisfied by any  control which  yields an optimal cost or pay-off. It consists of a system of forward-backward differential equations (e.g., state and adjoint equations)  and the extremum of a  Hamiltonian functional.     Since then the theory has been extensively developed in many directions, such as,  optimal control theory for finite  dimensional deterministic systems   with regular controls    \cite{cesari1983,ahmed2006}, where one can find a  broad and deep  generalization   of the classical Pontryagin  minimum  (equivalently maximum) principle   for deterministic systems.  The minimum principle is also extended to infinite dimensional systems, see \cite{ahmed1981,Fattorini1999,ahmed2005,ahmed-charalambous2007} and the references cited therein.\\

The stochastic minimum principle is another important extension of the Pontryagin minimum principle for systems subject to probabilistic randomness. In the stochastic case, there are basically different approaches based on the assumptions employed to derive the stochastic minimum principle. Specifically,  \cite{kushner1972} utilizes spike variations and Neustadt's variational principle,   \cite{haussmann1986}
utilizes Girsanov's measure transformation for non degenerate controlled diffusion processes, while \cite{bismut1978} utilizes the martingale representation to derive the adjoint equation. The martingale representation approach is further developed in \cite{bensoussan1981,bensoussan1983}. Further results utilizing the martingale representation approach are established  in \cite{peng1990} for control dependent diffusion processes utilizing second-order variations leading to a minimum principle which differs from the deterministic case in the sense that the effect of control dependent diffusion terms are fully explored. Subsequent extensions are given in \cite{cadenillas-karatzas1995} for stochastic systems with random coefficients, in \cite{elliott-kohlmann1994} utilizing stochastic flows to derive results similar to
\cite{peng1990}, and in  \cite{zhou1991} establishing  relationships between stochastic minimum principle and dynamic programming. The martingale approach to stochastic minimum principle sparked the interest in studying backward and forward stochastic differential equations. An excellent account on the stochastic minimum principle is found in \cite{yong-zhou1999} which also includes an anthology of references. Extensions of the stochastic maximum principle for relaxed controls using the topology of weak convergence  are found in \cite{mezerdi-bahlali2002,bahlali-mezerdi-djehiche2006,bahlali-djehiche-mezerdi2007a,bahlali-djehiche-mezerdi2007b,bahlali2008}, where relations to strict controls are also investigated.  \\
Recent developments and extensions are found in  \cite{ buckdahn-djehiche-li2011,zhang-shi2011,zhang-elliott-siu2012}  and references therein. \\ The area of mathematical finance, specifically portfolio optimization,  has utilized the stochastic minimum principle extensively to derive optimal strategies.  \\

In general, the stochastic minimum principle is specific to the information structures available to the control. Specifically,  in applications of control theory, there are many problems in  physical sciences  and engineering,   where systems are modeled by stochastic differential equations driven by controls which are also  stochastic processes with specific information structure, such as, full information or partial information.  Mathematically, information structures are modeled via the minimal sigma algebra generated by the available information process, and it is this process that the controller uses to generate control actions.    For full-information problems in which the information structure is Markovian,   one  often employs  Bellman's  principle of optimality to  construct, what is known as, the  HJB (Hamilton-Jacob-Bellman) equation, a nonlinear PDE defined on the state space of the system under investigation.  This equation describes the evolution of the   value function which is used to  construct the state feedback control law   provided this  function  is at least once differentiable with respect to the state variable.    This however requires solving the HJB equation   which may have a viscosity solution but  not  sufficiently  smooth \cite{yong-zhou1999}.  For non-Markovian controlled  diffusion systems with general information structures the HJB equation does not apply. For information structures which correspond  to full information or partial information the stochastic minimum principle is often employed \cite{bensoussan1983,charalambous-hibey1996,ahmed-charalambous2007}, although the partial information case is mathematically more demanding. However, this line of research is feasible provided existence of optimal controls is guaranteed. For non convex  control problems, it is well  known that the problem may have no optimal solution if the admissible controls are merely measurable functions with values in the set  $U$ which is non convex.   Nevertheless, this problem can be partially overcome by introducing the  relaxed controls and then approximating the relaxed controls by the standard regular controls.\\

In this article we  consider  stochastic control systems with information structures corresponding to full information and partial information,  which are driven by relaxed controls.  Specifically,  controls which are conditional  probability distributions, measurable with respect to full or partial information.  We treat stochastic differential equations driven by both Brownian motion and Le\'vy process or Poisson jump process. We show   existence of optimal policies  among the class of  relaxed controls under general conditions, with respect to an  established  topology of weak$^*$ convergence.  Then we proceed with the derivation of stochastic minimum principle, for both  the full information and the partial information cases. The Hamiltonian system of equations is derived in a systematic manner utilizing   the semi martingale representation theorem and the Riesz representation theorem, leading very naturally to the existence of the adjoint processes satisfying a Backward stochastic differential equation in an appropriate space.  We also discuss the realizability of relaxed controls by regular controls using the Krein-Millman theorem. The methodology we consider is applied to stochastic differential equations driven by both Brownian motion and Poisson jump process.   The basic procedure follows the one introduced   in \cite{ahmed2005,ahmed2006} for deterministic systems, augmented by the martingale representation approach to stochastic control. The material presented for full information  compliment the previous work on relaxed controls   found in \cite{mezerdi-bahlali2002,bahlali-mezerdi-djehiche2006,bahlali-djehiche-mezerdi2007a,bahlali-djehiche-mezerdi2007b,bahlali2008}, where the authors utilize alternative methods to derive related results. \\

\par The rest of the paper is organized as follows. In section 2 we present some typical notations and formulate the optimal control problem considered in this paper. In section 3, we consider the question of  existence of optimal relaxed  controls.  Section 4  contains an  interesting fundamental  result  characterizing semi martingales. Here we construct a Hilbert space characterizing the space of semi martingales (starting from zero). This  is used later in the development of necessary conditions.  Section 5 is devoted to the development of necessary conditions of optimality. In section 6 we extend the previous results to cover stochastic systems driven by jump processes. In section 7, we specialize to regular controls and obtain the usual necessary conditions of optimality. In section 8 we address the question of realizability of relaxed controls by regular controls.    The  paper is concluded with some comments on possible extensions of our results.

 \section{ Formulation of  Stochastic Relaxed Control Problem}
 \label{formulation}

 In this  section we  introduce  the mathematical  model for the  stochastic control system and the pay-off functional as a measure of performance. The distinction between full and partial information structures are also presented. \\

  Let $(\Omega,{\cal F}, {\cal F}_{t\geq 0}, P)$ denote a complete filtered probability space where $\{ {\cal F}_t, t \geq 0\}$ is an increasing family of subsigma algebras of the $\sigma$-algebra ${\cal F}.$  For any random variable $z$, ${\cal E}(z) \equiv \int_{\Omega} z(\omega) P(d\omega)$ denotes the expected value (average) of the  random variable $z.$  Let $\{ W(t), t \geq 0\}$ denote the $R^m$-valued  standard Brownian motion with $P\{ W(0) = 0\} = 1$ defined on the filtered probability space   $(\Omega,{\cal F}, {\cal F}_{t\geq 0}, P)$. Let ${\cal G}_t \subset {\cal F}_t$ denote a family of sub-sigma algebras of the $\sigma$-algebra ${\cal F}_t, t \geq 0.$

  Let $I = [0,T]$ be any finite interval,  $U$  any closed bounded subset of $R^d$ and   ${\cal M}(U)$ the space of regular bounded  signed Borel measures   on ${\cal B}(U)$, the Borel subsets of $U$ and ${\cal M}_1(U) \subset {\cal M}(U)$ the space of regular probability measures.  Controls based on  partial information (respectively full information) will be described  through the topological dual of the Banach space   $L_1^a(I,C(U))$, the $L_1$-space of ${\cal G}_t$ (respectively ${\cal F}_t$) adapted $C(U)$ valued functions.  The  dual of this space is given by  $L_{\infty}^a(I,{\cal M}(U))$ which, for partial information, consists of weak star  measurable  ${\cal G}_t$ adapted ${\cal M}(U)$ valued functions (signed measures), while for full information it  consists of ${\cal F}_t$ adapted functions defined similarly.  For controls based on partial information  (respectively full information) we are interested in the subspace $L_{\infty}^a(I,{\cal M}_1(U)) \subset L_{\infty}^a(I,{\cal M}(U))$  of probability measure valued  ${\cal G}_t$ (respectively ${\cal F}_t$) adapted functions.  Let ${\cal U}_{ad} \equiv L_{\infty}^a(I,{\cal M}_1(U))$ denote the class of admissible controls, called the relaxed controls, where the distinction between full information and partial information is only specified in terms of the $\sigma$-algebras  ${\cal F}_t$ and ${\cal G}_t$, respectively. \par  We consider the following stochastic system in $R^n$ governed by the Ito differential equation which is  driven by relaxed control,   \begin{eqnarray*} dx(t) = \bigl(\int_{U}  b(t,x(t),\xi) u_t(d\xi)\bigr)~ dt + \bigl(\int_U \sigma(t,x(t),\xi)u_t(d\xi)\bigr)~ dW(t), x(0) = x_0, t \in I, \end{eqnarray*}  where  $b: I\times R^n\times U \longrightarrow R^n$ denotes the drift and $\sigma : I\times R^n\times U \longrightarrow {\cal L}(R^m,R^n)$ the diffusion parameters. For simplicity of notation we prefer to write the above controlled dynamic system in the form \begin{eqnarray} dx(t) =  b(t,x(t),u_t)~ dt + \sigma(t,x(t),u_t)~ dW(t), x(0) = x_0, t \in I, \label{eq1}\end{eqnarray} for any  $u \in {\cal U}_{ad}.$ The cost functional is given by \begin{eqnarray*} J(u) \equiv {\cal E}\biggl\{\int_0^T \ell(t,x(t),u_t) dt + \Phi(x(T))\biggr\}. \end{eqnarray*} The problem is to find a control $u^o \in {\cal U}_{ad}$ such that $J(u^o) \leq J(u)$ for all $u \in {\cal U}_{ad}.$  We consider the question of existence of optimal controls and characterization of such controls in the form of necessary conditions of optimality (Pontryagin minimum principle). For necessary conditions of optimality we follow the procedure  developed in \cite{ahmed2006}, pp.271-293 which we extend from  deterministic to stochastic  systems.

  \section{ Existence of Optimal Relaxed Controls}
  \label{existence}

\vskip6pt Consider the system (\ref{eq1}) with   $b$ and $\sigma$ denoting  the infinitesimal generators representing the drift and diffusion given  by the   Borel measurable  maps: \begin{eqnarray*} b:I\times R^n \times R^d \longrightarrow R^n , \sigma: I\times R^n \times R^d \longrightarrow {\cal L}(R^m,R^n).\end{eqnarray*}
We assume that they satisfy  the following basic properties:
 there exists a $K \in L_2^+(I)$(nonnegative functions belonging to  $L_2(I)$),  such that
  \begin{eqnarray*}  &~& \hbox{(A1)}: |b(t,x,\xi)-b(t,y,\xi)|_{R^n} \leq K(t) |x-y|_{R^n} ~\hbox{uniformly in}~ \xi \in U \\  &~& \hbox{(A2)}: |b(t,x,\xi)|_{R^n} \leq K(t) (1 + |x|_{R^n}) ~\hbox{uniformly in}~ \xi \in U\\  &~& \hbox{(A3)}: |\sigma(t,x,\xi)-\sigma(t,y,\xi)|_{{\cal L}(R^m,R^n)} \leq K(t) |x-y|_{R^n} ~\hbox{uniformly in}~ \xi \in U  \\  &~& \hbox{(A4)}:  |\sigma(t,x,\xi)|_{{\cal L}(R^m,R^n)} \leq K(t) (1+ |x|_{R^n}) ~\hbox{uniformly in}~ \xi \in U \\  &~& \hbox{(A5)}:  b(t,x,\cdot), \sigma(t,x,\cdot) ~\hbox{are continuous in}~ \xi \in U   ~\hbox{uniformly in} \: t\in [0,T], x\in R^n.\end{eqnarray*}

   For admissible controls, we choose the set of relaxed controls given by  ${\cal U}_{ad} \equiv L_{\infty}^a(I,{\cal M}_1(U))$  which are  stochastic processes, adapted to a given sigma algebra (to be specified later), and taking values in the space of probability measures $M_1(U).$  This is endowed with the weak star topology also called vague topology.  A sequence $u^n \in {\cal U}_{ad}$ is said to converge  vaguely  to $u^o,$ written $ u^n\buildrel v\over\longrightarrow u^o$, iff for every  $\varphi \in L_{ 1}^a(I,C(U))$
  $$ {\cal E} \int_{I\times U}  \varphi(t,\xi) u_t^n(d\xi) dt \rightarrow  {\cal E} \int_{I\times U}  \varphi(t,\xi) u_t^o(d\xi) dt ~~ \hbox{as}~~ n \rightarrow \infty.$$ With respect to this vague (weak star)  topology, ${\cal U}_{ad}$ is  compact and from here on we assume that ${\cal U}_{ad}$ has been endowed with this vague topology.

     \vskip6pt Let $B_{\infty}^a(I, L_2(\Omega,R^n))$ denote the space of ${\cal F}_t$-adapted $R^n$ valued second order random processes endowed with the norm topology  $\parallel  \cdot \parallel$ given by  $$ \parallel x\parallel^2  \equiv \sup \{ {\cal E}|x(t)|_{R^n}^2, t \in I\}.$$
\vskip6pt  With this preparation, we can now  present the following lemma proving    existence of solutions and their continuous dependence on controls.

\vskip6pt\noindent {\bf Lemma 3.1} Consider the controlled stochastic differential equation (\ref{eq1}) and suppose the assumptions (A1)-(A5) hold. Then for any ${\cal F}_0$-measurable initial state $x_0$ having finite second moment, and any $u \in {\cal U}_{ad}$,  the system (\ref{eq1}) has a unique solution   $x \in B_{\infty}^a(I,L_2(\Omega,R^n))$  having continuous modification. In other words, $x \in C(I,R^n)$ P-a.s. Further, the solution is continuously dependent on the control in the sense that as $u^n \buildrel v \over\longrightarrow u^o$  in ${\cal U}_{ad}$, the corresponding solutions   $x^n \buildrel s \over\longrightarrow x^o $ in $B_{\infty}^a(I,L_2(\Omega,R^n).$

\vskip6pt \noindent {\bf Proof.} The proof for the first part of the lemma is classical  and hence we present only an outline. It is based on the  Banach fixed point theorem  applied to the operator  $F$  on the Banach space  $ B_{\infty}^a(I,L_2(\Omega,R^n))$ where
\begin{eqnarray} (Fx)(t) \equiv x_0 + \int_0^t b(s,x(s),u_s)~ds + \int_0^t \sigma(s,x(s),u_s)~dW(s), t \in I\equiv [0,T]. \label{eq2} \end{eqnarray} Under the assumptions (A1)-(A4), it is easy to verify using classical martingale inequality that  $F: B_{\infty}^a(I,L_2(\Omega,R^n))$ to itself. Then using the metric $d$ given by $d = d_T$ where  $$ d_t^2 (x,y) \equiv \sup\{ {\cal E} |x(s)-y(s)|_{R^n}^2, 0 \leq s \leq t\}$$ for  $t \in I,$ one can verify that the  $n-th$ iterate of $F$ denoted by $F^n \equiv FoF\cdots oF$ ($n$ times) is a contraction. Then by Banach fixed point theorem $F^n$ has a unique fixed point in $B_{\infty}^a(I,L_2(\Omega,R^n))$ and hence  $F$ itself has one and the same fixed point [1]. The continuity of the sample paths however follows from classical Borel-Canteli lemma.  Now consider the second part asserting the  continuity of the control to solution map  $u\longrightarrow x.$ For this  one proceeds as follows. Suppose the assumption  (A5) holds and   let $\{u^n,u^o\}$ be any sequence  of  controls from ${\cal U}_{ad}$ and  $\{x^n,x^o\}$ denote the corresponding sequence  of solutions of the system (\ref{eq1}).  Let $u^n \buildrel v\over\longrightarrow u^o.$ We must show that  $x^n \buildrel s \over\longrightarrow x^o$ in  $B_{\infty}^a(I, L_2(\Omega,R^n)).$  We present only a hint. Using the definition of solution, it is easy to verify that \begin{eqnarray}&~&  x^n(t)-x^o(t) = \int_0^t [ b(s,x^n(s),u^n_s)-b(s,x^o(s),u_s^n)] ds  \nonumber \\ &~& ~~~~~~~~~~+ \int_0^t [\sigma(s,x^n(s),u^n_s)-\sigma(s,x^o(s),u_s^n)] dW(s) + e_{1,n}(t) + e_{2,n}(t), t \in I \nonumber \\ \label{eq3}.\end{eqnarray} where
\begin{eqnarray*}&~&  e_{1,n}(t) =  \int_0^t [ b(s,x^o(s),u^n_s)-b(s,x^o(s),u_s^o)] ds \\ &~&  e_{2,n}(t) = \int_0^t [ \sigma(s,x^o(s),u^n_s)-\sigma(s,x^o(s),u_s^o)] dW(s).\end{eqnarray*} Using the standard martingale inequality it follows from this that there exist constants $C_1,C_2> 0$ such that
\begin{eqnarray} E|x^n(t)-x^o(t)|^2 \leq C_1 \int_0^t K^2(s)  E|x^n(s)-x^o(s)|^2 + C_2 \bigl( {\cal E}|e_{1,n}|^2 + {\cal E}|e_{2,n}|^2\bigr). \label{eq4} \end{eqnarray} Clearly, $$ {\cal E}|e_{1,n}|^2 \leq T {\cal E}\int_0^t  |b(s,x^o(s),u_s^n)-b(s,x^o(s),u_s^o)|_{R^n}^2 ds $$ and  $$ {\cal E}|e_{2,n}|^2 \leq ~ 4 {\cal E}\int_0^T |\sigma(s,x^o(s),u_s^n)-\sigma(s,x^o(s),u_s^o)|_{{\cal L}(R^m,R^n)}^2 ds. $$ Now note that by virtue of  vague  convergence of $u^n$ to $u^o$, the  integrands of  the above inequalities converge to zero  for almost all $s \in I,$ P-a.s  and it follows from (A2) and (A4) that they are dominated by integrable functions. So by Lebesgue dominated convergence theorem  the integrals $\{e_{1,n},e_{2,n}\}$ converge to zero uniformly on $I.$  The assertion then follows  from  Gronwall inequality, applied to the inequality (\ref{eq4}). This completes the outline.  $\bullet$

\vskip6pt\noindent{\bf Optimal Control Problem.} Consider the controlled system (\ref{eq1}) and the cost functional given by \begin{eqnarray} J(u) \equiv {\cal E} \bigl\{ \int_0^T \ell(t,x(t),u_t) dt + \Phi(x(T))\bigr\} \label{eq5} \end{eqnarray}  where $\ell$ and $\Phi$ are suitable functions which are measures of mismatch between the desired flow and the flow that results from the choice of the control $u.$ The problem, as stated in section 2,  is to find a control from the class of admissible (relaxed) controls ${\cal U}_{ad}$ that minimizes the functional (\ref{eq5}). We present the following existence result.

\vskip6pt\noindent {\bf Theorem 3.2} Consider the control problem as stated above. Suppose the assumptions of Lemma 3.1 hold, and further suppose  $\ell: I\times R^n\times U \longrightarrow (-\infty,+\infty]$ and $\Phi:R^n \longrightarrow (-\infty,+\infty]$ are Borel measurable maps satisfying the following conditions: \par (a1): $x\longrightarrow \ell(t,x,\xi)$ is  continuous on $R^n$ for each $t\in I$, uniformly with respect to $\xi \in U.$ \par (a2): $\exists$ $h \in L_1^+(I)$ such that $|\ell(t,x,\xi)| \leq h(t) (1 + |x|_{R^n}^2)$\par(a3): $x \longrightarrow \Phi(x)$ is lower semi continuous on $R^n$ and  $\exists$ $c_0,c_1\geq 0$ such that $|\Phi(x)| \leq c_0 + c_1 |x|_{R^n}^2.$ \\ Then,  there exists an optimal control $u \in {\cal U}_{ad}$ at which $J$ attains its minimum.\vskip6pt \noindent {\bf Proof.} Since ${\cal U}_{ad}$ is compact in the vague  topology, it suffices to prove that $J$ is lower semi continuous with respect to this topology.  Suppose $u^n \buildrel v \over\longrightarrow u^o$ in ${\cal U}_{ad}$  and let $\{x^n,x^o\} \subset B_{\infty}^a(I,L_2(\Omega,R^n))$ denote the solutions of equation (\ref{eq1}) corresponding to the sequence of controls $\{u^n,u^o\} \subset {\cal U}_{ad}$. Then by Lemma 3.1, along a subsequence if necessary, $x^n \buildrel s \over\longrightarrow x^o$ in $B_{\infty}(I,L_2(\Omega,R^n)).$ First note  that,  in view of the strong convergence, along  a subsequence if necessary, $x^n(T) \rightarrow x^o(T)$ P-a.s. Thus it follows from assumption  (a3) and Fatou's Lemma  that \begin{eqnarray} {\cal E} \{\Phi(x^o(T))\} \leq \liminf_n {\cal E} \{\Phi(x^n(T))\}. \label{eq6}\end{eqnarray} Considering  the running cost, it is easy to see that \begin{eqnarray} &~& \hskip-50pt {\cal E} \int_{I} \ell(t,x^o(t),u_t^o)~dt =  {\cal E}  \int_{I}\ell(t,x^o(t),u_t^o-u_t^n)~dt \nonumber \\ &~&  +  {\cal E} \int_{I} ( \ell(t,x^o(t), u_t^n)-\ell(t,x^n(t),u_t^n))~dt +  {\cal E} \int_I \ell(t,x^n(t),u_t^n)~dt.\label{eq7} \end{eqnarray} By virtue of vague  convergence of $u^n$ to $u^o,$  it is evident that  for every $\varepsilon >0$ there exists an integer $n_{1,\varepsilon}$ sufficiently large, such that the absolute value of the  first term on the right hand side of equation (\ref{eq7})  is less than $\varepsilon/2 $ for all $n \geq n_{1,\varepsilon}.$ By virtue of assumption (a1)-(a2), in  particular the continuity of $\ell$ in $x$ uniformly in $U$,  it is easy to verify that there exists an integer $n_{2,\varepsilon}$ such that for all $n \geq n_{2,\varepsilon}$, the absolute value of the second term on the right hand side is less than $\varepsilon/2.$  By combining these facts we obtain the following inequality
\begin{eqnarray*} &~& \hskip-50pt {\cal E} \int_{I} \ell(t,x^o(t),u_t^o)~dt  \leq \varepsilon + \int_I \ell(t,x^n(t),u^n_t) dt  \end{eqnarray*} for all $n \geq n_{1,\varepsilon} \bigvee n_{2,\varepsilon}.$ Since $\varepsilon >0$ is otherwise arbitrary, it follows from the above inequality that \begin{eqnarray} &~& \hskip-50pt {\cal E} \int_{I} \ell(t,x^o(t),u_t^o)~dt \leq \liminf_{n}  {\cal E} \int_I \ell(t,x^n(t),u_t^n)~dt. \label{eq8} \end{eqnarray} Combining (\ref{eq6}) and (\ref{eq8}) we arrive at the conclusion that $J(u^o) \leq \liminf_{n} J(u^n)$ thereby proving  lower semi continuity of $J$ in the vague topology.  Since ${\cal U}_{ad}$ is  compact in this vague  topology, $J$ attains its minimum on it.  This proves the existence of an optimal control. $\bullet$

\  \
Note that the  existence is  proved  under general conditions, irrespectively of  whether the information structure to the control is full  or partial.

\section{ Construction of a Hilbert  Space of Semi Martingales}
\label{martingales}

In the preceding section we have presented a result on  existence of optimal controls.   In the following  section  we consider the problem of characterizing optimal controls in the form of necessary conditions of optimality. For  this we shall utilize martingale approach hence we need to consider certain fundamental   properties of semi martingales.  These properties are studied in  this section.  Before we consider such properties, we wish to provide the  technical reasons for their study.  Consider the system (\ref{eq1}) with the cost functional (\ref{eq5}) and the admissible controls ${\cal U}_{ad}\equiv L_{\infty}^a(I,{\cal M}_1(U)) $ as described above. Recall that these are either ${\cal F}_t$ or ${\cal G}_t$-adapted probability measure valued random processes, depending on whether the information structure used to construct the controls is full or partial.    For the  necessary conditions of optimality  we need stronger regularity properties for the drift and diffusion parameters $\{ b,\sigma\}$ as well as the cost integrands $\{\ell,\Phi\}.$  They are presented as follows: \vskip6pt\noindent  {\bf (NC1):}  The triple $\{b,\sigma,\ell\} $ are measurable in $t \in I$, and the quadruple   $\{b,\sigma,\ell, \Phi\} $ are  once continuously differentiable with respect to the state variable $x \in R^n.$  The first spatial derivatives of $\{b,\sigma\}$ are bounded uniformly on $I\times R^n \times U.$ \vskip6pt Considering the Gateaux derivative of $\sigma$ with respect to the state variable  at the point $(t,z,\nu) \in I\times R^n\times M_1(U)$   in the direction $\eta \in R^n$ we have
$$ \lim_{\varepsilon \rightarrow 0} (1/\varepsilon) ( \sigma(t,z + \varepsilon \eta, \nu)- \sigma(t,z,\nu)) \equiv \sigma_x(t,z,\nu; \eta).$$ Note that $\eta \longrightarrow \sigma_x(t,z,\nu; \eta)$ is linear and it follows from the assumption (NC1) that  there exists a finite positive number $\beta$ such that
$$ |\sigma_x(t,z,\nu; \eta)|_{{\cal L}(R^m,R^n)} \leq \beta |\eta|_{R^n}.$$
In order to present the necessary conditions of optimality we need the so called variational equation.
Suppose $u^o \in {\cal U}_{ad}$ denote the optimal control and $u \in {\cal U}_{ad}$ any other control.  Since ${\cal U}_{ad}$ is convex, for any $\varepsilon \in [0,1]$, the control $$ u^{\varepsilon} \equiv u^o + \varepsilon (u-u^o) \in {\cal U}_{ad}.$$ Let $x^{\varepsilon}, x^{o} \in B_{\infty}^a(I,L_2(\Omega,R^n)) $ denote the solutions of the system equation (\ref{eq1}) corresponding to the controls $u^{\varepsilon}$ and $u^o$ respectively.  Consider the limit $$ y \equiv \lim_{\varepsilon\downarrow 0} (1/\varepsilon) (x^{\varepsilon}-x^o). $$
We have the following result characterizing the process $y.$

\vskip6pt\noindent {\bf Lemma 4.1}
The process $y$ is an element of the Banach space    $B_{\infty}^a(I,L_2(\Omega,R^n))$ and it  is the unique solution of the variational SDE \begin{eqnarray} &~& dy(t) = b_x(t,x^o(t),u_t^o)~y(t)~ dt + \sigma_x(t,x^o(t),u_t^o; y(t))~dW (t)\nonumber \\ &~& ~~~~~~~~~~~~~~~~~~~~~~~~~~~~~~ + b(t,x^o(t),u_t-u_t^o)~dt + \sigma(t,x^o(t),u_t-u_t^o)~dW(t),\label{eq9}   \\ &~& y(0) = 0\nonumber, \end{eqnarray} having a continuous modification.

\vskip6pt\noindent {\bf Proof.}  This is a linear SDE and so one can have  a  closed form solution. Indeed, considering the homogenous part given by  $$dz(t) = b_x(t,x^o(t),u_t^o)~z(t)~ dt + \sigma_x(t,x^o(t),u_t^o; z(t))~dW(t), z(s) = \zeta , 0 \leq s\leq t < \infty,$$ it follows from the assumption (NC1) and  Lemma 3.1 that it has a unique solution $z$ given by $$ z(t) = \Psi(t,s) \zeta, t \geq s, $$ where $\Psi(t,s), 0 \leq s \leq t<\infty  $ is the random (${\cal F}_t$ measurable)  transition operator for the homogenous system. Since the spatial derivatives of $b$ and $\sigma$ are uniformly bounded, the transition operator $\Psi(t,s), 0 \leq s \leq t \leq T $  is uniformly  $P$ almost surely bounded (with values in the space of $n\times n$ matrices).  Considering the non homogenous system (\ref{eq9}),  the solution is then  given by  \begin{eqnarray}  y(t) = \int_0^t \Psi(t,s)d\eta(s)  \label{eq10} \end{eqnarray} where $\eta$ is the semi martingale given by
\begin{eqnarray}  d\eta(t) = b(t,x^o(t),u_t-u_t^o)~dt +  \sigma(t,x^o(t),u_t-u_t^o)~dW(t), \eta(0) = 0.\label{eq11}\end{eqnarray} Note that $\eta$ is a continuous  square integrable   ${\cal F}_t$  semi martingale. This proves  the existence, uniqueness and regularity property of the solutions of system (\ref{eq9}).  This is one approach. An alternate   approach is the same as that of Lemma 3.1.  Here one notes that the drift and the diffusion terms of equation (\ref{eq9}) satisfy the basic assumptions of Lemma 3.1.  So the existence of a solution follows from the Banach fixed point theorem as in lemma 3.1.  The fact that it has continuous modification  follows directly from the  representation (\ref{eq10}) and the continuity of the semi martingale $\eta.$  $\bullet$

\vskip6pt Later in the sequel we need certain important and  interesting properties  of semi martingales. Let $L_2^a(I,R^n) \subset L_2(I\times \Omega,R^n)$ denote the space of ${\cal F}_t$-adapted random processes $\{v(t), t \in I\}$   such that $$ {\cal E}\int_{I} |v(t)|_{R^n}^2 dt < \infty.$$  Similarly, let $L_2^a(I,{\cal L}(R^m,R^n)) \subset L_2(I\times \Omega ,{\cal L}(R^m,R^n))$ denote the space of  ${\cal F}_t$-adapted $n\times m$ matrix valued random processes $\{ \Sigma(t), t \in I\}$ such that
 $$ {\cal E}\int_{I} |\Sigma(t)|_{{\cal L}(R^m,R^n)}^2 dt = {\cal E} \int_I tr(\Sigma^*(t)\Sigma(t)) dt < \infty.$$  Since $I$ is a finite interval , it is clear that  $B_{\infty}^a(I,L_2(\Omega,R^n)) \subset L_2^a(I,R^n).$

\vskip6pt\noindent{\bf Definition 4.2}  An $R^n$-valued  random process $\{m(t), t \in I\}$ is said to be a square integrable continuous   ${\cal F}_t$-semi martingale iff it is representable in the form
\begin{eqnarray} m(t) = m(0) + \int_0^t v(s) ds + \int_0^t \Sigma(s) dW(s), t \in I, \label{eq12} \end{eqnarray}  for some $v \in L_2^a(I,R^n)$ and $\Sigma \in L_2^a(I,{\cal L}(R^m,R^n))$  and for some $R^n$-valued ${\cal F}_0$ measurable  random variable  $m(0)$ having finite second moment.

 \vskip6pt  We  introduce the following class of ${\cal F}_t$-semi martingales:  \begin{eqnarray} &~& \hskip-40pt  {\cal SM}_0^2 \equiv \biggl \{  m: m(t)  = \int_0^t v(s) ds + \int_0^t \Sigma(s) dW(s),  t \in I, \nonumber \\ &~& ~~~~~~~~~~~~~~~~~~~~for ~  v \in L_2^a(I,R^n)~\hbox{ and}~ \Sigma \in L_2^a(I,{\cal L}(R^m,R^n)) \biggr \}.\label{eq13} \end{eqnarray} \vskip6pt  Now we present a fundamental result which has the potential of many other applications.

\vskip6pt\noindent{\bf Theorem  4.3} The class ${\cal SM}_0^2$ is a real linear vector space and it is a  Hilbert space with respect to the norm topology   $\parallel m\parallel_{{SM}^2_0}$ arising from
$$ \parallel m \parallel_{{\cal SM}^2_0}^2  \equiv {\cal E}\int_I |v(t)|_{R^n}^2 dt + {\cal E} \int_{I} tr(\Sigma^*(t)\Sigma(t)) dt. $$ Further, the space ${\cal SM}_0^2$ is isometrically isomorphic to $L_2^a(I,R^n)\times L_2^a(I,{\cal L}(R^m,R^n)),$ written as ${\cal SM}_0^2 \cong L_2^a(I,R^n)\times L_2^a(I,{\cal L}(R^m,R^n)).$

\vskip6pt\noindent {\bf Proof}  Note that each $m \in {\cal SM}_0^2$ corresponds to a  pair  $$(v,\Sigma)\in L_2^a(I,R^n)\times L_2^a(I,{\cal L}(R^m,R^n)).$$   We may call the pair $(v,\Sigma)$ the infinitesimal generator (or simply the intensity) of the semi martingale  $m.$   Let $m_1 \in {\cal SM}_0^2$ corresponding to the intensity process  $(v_1,\Sigma_1) \in  L_2^a(I,R^n)\times L_2^a(I,{\cal L}(R^m,R^n))$ and $m_2 \in {\cal SM}_0^2$  corresponding to the intensity process  $(v_2,\Sigma_2) \in  L_2^a(I,R^n)\times L_2^a(I,{\cal L}(R^m,R^n))$ respectively. Clearly, $v_1+v_2 \in L_2^a(I,R^n)$ and $\Sigma_1+\Sigma_2 \in L_2^a(I,{\cal L}(R^m,R^n)).$ Hence $m \equiv  m_1+m_2,$ with intensity process $(v_1+v_2, \Sigma_1+ \Sigma_2),$ is an element of ${\cal SM}_0^2.$  For any  real number $\alpha$ and any $m \in {\cal SM}_0^2 $ with intensity process $(v,\Sigma) \in  L_2^a(I,R^n)\times L_2^a(I,{\cal L}(R^m,R^n)$ we have $\alpha m \in {\cal SM}_0^2$ with intensity process $(\alpha v,\alpha \Sigma) \in  L_2^a(I,R^n)\times L_2^a(I,{\cal L}(R^m,R^n).$ Thus ${\cal SM}_0^2$ is a linear vector space. We now furnish this with a scalar product and norm topology. Let $m_1,m_2 \in {\cal SM}_0^2$ with the intensity pairs $(v_1,\Sigma_1), (v_2,\Sigma_2)$ respectively  and define \begin{eqnarray} (m_1,m_2)_{{\cal SM}_0^2} \equiv {\cal E}\int_I(v_1(t),v_2(t)) dt + {\cal E} \int_I tr(\Sigma_1^*(t)\Sigma_2(t)) dt.\label{eq14} \end{eqnarray}  The reader can easily verify that this gives a scalar product.  Clearly taking $m_2 = m_1$ we have the norm square of $m_1$ given by
\begin{eqnarray} \parallel m_1 \parallel^2_{{\cal SM}_0^2} =  (m_1,m_1)_{{\cal SM}_0^2} \equiv {\cal E}\int_I|v_1(t)|_{R^n}^2 dt + {\cal E} \int_I |\Sigma_1(t)|_{{\cal L}(R^m,R^n)}^2 dt. \label{eq15} \end{eqnarray}  It is easy to verify that the above expression defines a norm (modulo the null space).  Thus ${\cal SM}_0^2$ is a scalar product space. To show that it is a Hilbert space, it suffices to verify that it is complete. Let $\{m_n\} \subset {\cal SM}_0^2$ be a  Cauchy sequence corresponding to the sequence of intensity pairs $\{ (v_n,\Sigma_n)\} \subset L_2^a(I,R^n)\times L_2^a(I,{\cal L}(R^m,R^n)).$ Let $p \geq 1$ and consider the expression  \begin{eqnarray} &~& \parallel m_{n+p} - m_n\parallel_{{\cal SM}_0^2} =  \biggl ({\cal E}\int_I|v_{n+p}(t)-v_n(t)|_{R^n}^2 dt + {\cal E} \int_I |\Sigma_{n+p}(t)-\Sigma_{n}|_{{\cal L}(R^m,R^n)}^2 dt \biggr)^{1/2}.\nonumber \end{eqnarray} Since $\{m_n\}$ is a Cauchy sequence, $\lim_{n \rightarrow \infty} \parallel m_{n+p} - m_n\parallel_{{\cal SM}_0^2}  = 0$  for every $p \geq 1 $ and hence $\{(v_n,\Sigma_n)\}$ is a  Cauchy sequence in $ L_2^a(I,R^n)\times L_2^a(I,{\cal L}(R^m,R^n))$. But the later spaces are Hilbert and hence there exists a unique  pair $(v_o,\Sigma_o) \in L_2^a(I,R^n)\times L_2^a(I,{\cal L}(R^m,R^n))$ to which $(v_n,\Sigma_n)$ converges in norm (along a subsequence if necessary). Define the process $m_0$ by $$  m_0(t) = \int_0^t v_0(s) ds + \int_0^t \Sigma_0(s) dW(s), t \in I.$$ Clearly this is a semi martingale belonging to ${\cal SM}^2_0$ and it is the unique limit of the sequence of semi martingales $\{m_n\}.$  This proves that  ${\cal SM}^2_0$ is complete and hence a Hilbert space. Now we claim that for every $m \in {\cal SM}_0^2$ there exists a unique pair $(v,\Sigma)\in  L_2^a(I,R^n)\times L_2^a(I,{\cal L}(R^m,R^n))$ such that
$$  m(t) = \int_0^t v(s) ds + \int_0^t \Sigma(s) dW(s), t \in I.$$ Suppose this is false and there exists another pair $(v_1,\Sigma_1) \in L_2^a(I,R^n)\times L_2^a(I,{\cal L}(R^m,R^n))$ giving the same semi martingale $m.$  This means that $$ 0 = \int_0^t (v(s)-v_1(s)) ds + \int_0^t (\Sigma(s)-\Sigma_1(s)) dW(s), t \in I, $$ which is the same as $$ \int_0^t (v(s)-v_1(s)) ds  =  \int_0^t (\Sigma_1(s)-\Sigma(s)) dW(s), ~~ t \in I. $$ But this is impossible since a martingale can never equal a function of bounded variation. Hence $v_1 = v$ and $\Sigma_1 = \Sigma.$  Thus to every $m \in {\cal SM}_0^2$ there corresponds a unique pair $(v,\Sigma)\in L_2^a(I,R^n)\times L_2^a(I,{\cal L}(R^m,R^n))$ and conversely. The isometry follows from the expression (\ref{eq15}).  Hence ${\cal SM}_0^2 \cong L_2^a(I,R^n)\times L_2^a(I,{\cal L}(R^m,R^n)).$  This completes the proof. $ \bullet $

\section{Necessary Conditions of Optimality}
\label{necessary}
 Now we are prepared to develop the necessary conditions of optimality. The theory  of relaxed controls is found to be a powerful technique for developing necessary conditions of optimality for deterministic systems  \cite{ahmed2006},  Theorem 8.3.5. Here we use the same technique for systems governed by stochastic differential equations driven by  relaxed controls. \\
 Below, we provide the main theorem. Later we use this result to derive a simplified minimum principle for   both  full as well as  partial information.

 \vskip6pt\noindent {\bf Theorem 5.1} Consider the system (\ref{eq1}) and the cost functional (\ref{eq5}). An element $ u^o \in {\cal U}_{ad},$ with the corresponding solution $x^o \in B_{\infty}^a(I, L_2(\Omega,R^n))$ to be optimal, it is necessary that there exists a semi martingale  $m^o \in {\cal SM}_0^2$ with the intensity process $(\psi,Q) \in  L_2^a(I,R^n)\times L_2^a(I,{\cal L}(R^m,R^n))$ such that the following inequality and the equations (SDE) hold:

\begin{eqnarray} (1): \hskip-20pt &~& {\cal E} \int_0^T \bigl\{ (b(t,x^o(t),u_t),\psi(t)) + tr(Q^*(t)\sigma(t,x^o(t),u_t)) + \ell(t,x^o(t),u_t)\bigr\} dt \nonumber \\ &~& ~~~~~~~~~~\geq {\cal E} \int_0^T \bigl\{ (b(t,x^o(t),u_t^o),\psi(t)) + tr(Q^*(t)\sigma(t,x^o(t),u_t^o)) + \ell(t,x^o(t),u_t^o)\bigr\} dt  \nonumber  \\ &~& ~~\hbox{ for all }~~ u \in {\cal U}_{ad}. \label{eq16}\end{eqnarray}

\begin{eqnarray} (2): \hskip-20pt  &~&    dx^o(t) = b(t,x^o(t),u_t^o) dt + \sigma(t,x^o(t),u_t^o))dW(t) \nonumber \\ &~& x^o(0) = x_0 \label{eq17}  \end{eqnarray}

 \begin{eqnarray} (3):  \hskip-20pt &~&    -d \psi(t)  = b_x^*(t,x^o(t),u_t^o)\psi(t)  dt + V_{Q}(t) dt +\ell_x(t,x^o(t),u_t^o) dt - Q(t) dW(t)  \nonumber \\ &~& \psi(T) = \Phi_x(x^o(T)) \label{eq18}  \end{eqnarray} where  $V_{Q} \in L_2^a(I,R^n)$ is   given by  $(V_{Q}(t),\zeta) = tr (Q^*(t)\sigma_x(t,x^o(t),u_t^o; \zeta)), t \in I.$

\vskip6pt \noindent {\bf Proof} Suppose $u^o \in {\cal U}_{ad}$ is  the optimal control and $u \in {\cal U}_{ad}$ any other control.  Since ${\cal U}_{ad}$ is convex, for any $\varepsilon \in [0,1]$, the control $ u^{\varepsilon} \equiv u^o + \varepsilon (u-u^o) \in {\cal U}_{ad}.$ Let $x^{\varepsilon}, x^{o} \in B_{\infty}^a(I,L_2(\Omega,R^n)) $ denote the (strong)  solutions of the system equation (\ref{eq1}) corresponding to the controls $u^{\varepsilon}$ and $u^o$ respectively. Since $u^o$ is optimal it is clear that \begin{eqnarray}  J(u^{\varepsilon})-J(u^o) \geq 0 \label{eq19} \end{eqnarray}  for all $\varepsilon \in [0,1]$ and for all $u \in {\cal U}_{ad}.$  Let $dJ(u^o,u-u^0)$ denote the  Gateaux differential  of $J$ at $u^o$ in the direction $u-u^o.$ Dividing the expression (\ref{eq19}) by $\varepsilon$ and letting $\varepsilon \downarrow 0$ it is easy to verify that  \begin{eqnarray} dJ(u^o,u-u^0) = L(y) + {\cal E} \int_0^T \ell(t,x^o(t),u_t-u_t^o) dt \geq 0, ~~\forall~~ u \in {\cal U}_{ad} \label{eq20} \end{eqnarray}

where $L(y)$ is given by the functional  \begin{eqnarray} L(y) = {\cal E}\biggl\{  \int_0^T (\ell_x(t,x^o(t),u_t^o),y(t))~ dt + (\Phi_x(x^o(T)),y(T))\biggr\}. \label{eq21}\end{eqnarray} Since by Lemma 4.1, the process $y \in B_{\infty}^a(I,L_2(\Omega,R^n))$ and it is also   continuous P-a.s   it follows from  assumption (a2) of Theorem 3.2 and the assumption (NC1), that $y \longrightarrow L(y)$ is a continuous linear functional. Further, by Lemma 4.1, $\eta \longrightarrow y$ is a continuous linear map from the Hilbert space ${\cal SM}_0^2$ to the B-space $B_{\infty}^a(I,L_2(\Omega,R^n)) $  given by the expression (\ref{eq10}). Thus the composition map  $\eta \longrightarrow y \longrightarrow L(y) \equiv \tilde L (\eta)$ is a continuous linear functional on ${\cal SM}_0^2.$ Then by virtue of the classical Riesz representation theorem for Hilbert spaces, there exists a semi martingale $\varrho \in {\cal SM}_0^2$ with intensity $(\psi,Q) \in  L_2^a(I,R^n)\times L_2^a(I,{\cal L}(R^m,R^n))$ such that \begin{eqnarray} &~&  L(y) \equiv \tilde L (\eta) = (\varrho,\eta)_{{\cal SM}_0^2} =  {\cal E} \int_0^T (\psi(t), b(t,x^o(t),u_t-u_t^o)) dt \nonumber \\ &~& ~~~~~~~~~~~~~~~~~~~~~~~~~~~~~~~~~~~~~~~~~~+ {\cal E} \int_0^T tr (Q^*(t) \sigma(t,x^o(t),u_t-u_t^o)) dt . \label{eq22} \end{eqnarray} Substituting  the expression (\ref{eq22}) into the expression (\ref{eq20}) we obtain
\begin{eqnarray} &~&\hskip-50pt  dJ(u^o,u-u^0) = {\cal E} \int_0^T (\psi(t), b(t,x^o(t),u_t-u_t^o)) dt  \nonumber \\ &~&  ~~~~~~~~~~~+ {\cal E} \int_0^T tr (Q^*(t) \sigma(t,x^o(t),u_t-u_t^o)) \nonumber \\ &~& ~~~~~~~~~~~~~~~~~~~~~~~~~ + {\cal E} \int_0^T \ell(t,x^o(t),u_t-u_t^o) dt \geq 0, ~~\forall~~ u \in {\cal U}_{ad}. \label{eq23} \end{eqnarray}  The  necessary condition  given by the expression (\ref{eq16}) readily  follows from this.   Equation (\ref{eq17}) is the system equation along the optimal control state pair  $(u^o,x^o),$ so nothing to prove. We prove that the pair $ (\psi,Q)$ is given by the solution of the  adjoint equation   (\ref{eq18}). Computing  the Ito differential   of the scalar product $(y,\psi)$ we have the general expression \begin{eqnarray} d(y(t),\psi(t)) = (dy(t), \psi(t)) + (y(t),d\psi(t)) + <dy(t),d\psi(t)> \label{eq24}\end{eqnarray} where the last bracket  denotes the classical quadratic variation term.  Integrating this over $I = [0,T]$ and using the fact that  $y(0) = 0$, it follows from the variational equation (\ref{eq9}) that

\begin{eqnarray} &~& {\cal E} (y(T),\psi(T))=  {\cal E} \biggl\{ \int_0^T (y(t), b_x^*\psi(t) dt + \sigma_x^*(\psi(t))dW(t) + d \psi(t))  \nonumber \\ &~& ~~~ +\int_0^T (b^o,\psi(t)) dt + \int_0^T ((\sigma^o)^*\psi(t),dW(t))\biggr\}  + {\cal E} \int_0^T <dy(t),d \psi(t)>, \label{eq25}\end{eqnarray} where for convenience of notation we have used \begin{eqnarray*} &~& b_x \equiv b_x(t,x^o(t),u_t^o),~~ \sigma_x(\xi) \equiv  \sigma_x(t,x^o(t),u_t^o;\xi), \xi \in R^n,\\ &~& b^o = b(t,x^o(t),u_t-u_t^o),~~ \sigma^o \equiv \sigma(t,x^o(t),u_t-u_t^o).\end{eqnarray*} Note that the stochastic integrals in (\ref{eq25}) equal zero and hence make no contribution. This follows from the facts that  $\sigma_x^*(\psi(t)) \in L_2^a(I,{\cal L}(R^m,R^n))$ and $(\sigma^o)^*\psi \in L_2^a(I,R^m) $ as seen later.   So we can eliminate them giving  the following expression
\begin{eqnarray}  {\cal E} (y(T),\psi(T)) &=&  {\cal E} \biggl\{ \int_0^T (y(t), b_x^*\psi(t) dt + d \psi(t)) +  \int_0^T (b^o,\psi(t)) dt \biggr\} \nonumber  \\
& +& {\cal E} \int_0^T <dy(t),d \psi(t)>.  \label{eq26}\end{eqnarray}  Before we consider the quadratic variation term, let us recall that the Ito derivatives of the variation process $y$ and the adjoint process $\psi$ are of the following form:

 \begin{eqnarray*} dy(t)& =& \hbox{bounded variation terms} + \sigma_{x}(t,x^o(t),u_t^o; y(t)) dW(t) \nonumber \\
  &+& \sigma(t,x^o(t),u_t-u_t^o) dW(t) , \\ d \psi(t) &=& \hbox{bounded variation terms}  + Q(t) dW(t). \end{eqnarray*} Considering  now the quadratic variation term it is easy to verify that

\begin{eqnarray}  {\cal E} \int_0^T <dy(t),d \psi(t)> = {\cal E} \int_0^T \bigl \{  tr (Q^*(t)\sigma_x(y)) + tr(Q^*(t)\sigma^o) \bigr \} dt. \label{eq27} \end{eqnarray}  Clearly, the first term on the right hand side of the above expression is linear in $y.$  Thus there exists a process  $ V_Q(t), t \in I,$  given by the following expression  \begin{eqnarray}(V_Q(t),y(t)) \equiv tr(Q^*(t)\sigma_x(y)) \equiv tr ( Q^*(t) \sigma_x(t,x^o(t),u_t^o;y(t)) ). \label{eq28}  \end{eqnarray}  By assumption (NC1), $\sigma$ has uniformly bounded  spatial first derivative and it follows from the semi martingale representation Theorem 4.3  that $Q \in L_2^a(I,{\cal L}(R^m,R^n))$ and hence $V_{Q} \in L_2^a(I,R^n).$ Substituting (\ref{eq28}) into (\ref{eq27}) and then    (\ref{eq27}) into (\ref{eq26}), we obtain

\begin{eqnarray} &~& {\cal E} (y(T),\psi(T)) =  {\cal E} \biggl\{ \int_0^T (y(t), b_x^*\psi(t) dt + V_Q(t) dt - Q(t)dW(t) +  d \psi(t)) \nonumber \\ &~& ~~~~~~~~~~~~~~~~~~~~~~~~~~~~~~~~~~~~~~~~~~~~ +  \int_0^T (b^o,\psi(t)) dt  + tr(Q^*(t)\sigma^o) dt \biggr\}. \nonumber \\ &~&  \label{eq29}\end{eqnarray}
By setting \begin{eqnarray} &~&  b_x^*(t,x^o(t),u_t^o)\psi(t) dt +V_{Q}(t) dt  -Q(t) dW(t) + d\psi(t) = -\ell_x(t,x^o(t),u_t^o) dt \nonumber \\ &~& \psi(T) = \Phi_x(x^o(T)),  \label{eq30} \end{eqnarray}  it follows from (\ref{eq29}) and the expression for the functional $L$ given by (\ref{eq21})
that
\begin{eqnarray}  L(y) &=& {\cal E} (y(T),\psi(T)) + {\cal E}  \int_0^T (y(t),\ell_x(t,x^o(t),u_t^o)) dt \nonumber \\
& =&   {\cal E}  \int_0^T \{(b(t,x^o(t),u_t-u_t^o)\psi(t)) + tr(Q^*(t)\sigma(t,x^o(t),u_t-u_t^o))\} dt .   \label{eq31}\end{eqnarray} This is precisely what was obtained by the semi martingale argument giving (\ref{eq22}). Thus the pair $(\psi,Q)$ must satisfy the backward stochastic differential equation (\ref{eq30}) which is precisely the adjoint equation given by (\ref{eq18}) as  stated. Since $\psi$ satisfies the stochastic differential equation and $T$ is finite,  it follows from the classical  theory of Ito differential equations that $\psi$  is actually   an element of $B_{\infty}^a(I,L_2(\Omega,R^n)) \subset L_2^a(I,R^n).$ In other words, $\psi$ is more regular than predicted by semi martingale theory.  Hence by our assumption on $\sigma$ it is easy to verify that   $\sigma_x^*(\psi) \in L_2^a(I,{\cal L}(R^m,R^n))$ and $ (\sigma^o)^*\psi \in L_2^a(I,R^m) $ as stated before.   Thus we have completed the proof.  $\bullet$

\vskip6pt\noindent {\bf Remark 5.2} Define the Hamiltonian $$ H:I\times R^n\times R^n\times {\cal L}(R^m,R^n)\times {\cal  M}_1(U) \longrightarrow R $$ by  $$ H(t,\xi,\zeta,M,\nu) = (b(t,\xi,\nu),\zeta) + tr (M^*\sigma(t,\xi,\nu)) + \ell(t,\xi,\nu).$$  In terms of this Hamiltonian, the necessary conditions of optimality  (\ref{eq16})-(\ref{eq18}) can be written compactly as follows \begin{eqnarray} &~& {\cal E} \int_0^T H(t,x^o(t),\psi(t),Q(t),u_t) dt  \geq   {\cal E} \int_0^T H(t,x^o(t),\psi(t),Q(t),u_t^o) dt ~~ \nonumber \\ &~&  \hbox{for all}~ u \in {\cal U}_{ad},\label{eq32} \end{eqnarray} where the triple  $\{x^o,\psi,Q\}$ is  the unique solution of the following  Hamiltonian system
  \begin{eqnarray} &~&\hskip-20pt  dx^o(t) = H_\psi (t,x^o(t),\psi(t),Q(t),u_t^o)     dt + \sigma(t,x^o(t),u_t^o) dW(t),~~ x^o(0) = x_0, \label{eq33} \\ \nonumber \\ &~&\hskip-20pt  d \psi(t) = - H_x (t,x^o(t),\psi(t),Q(t),u_t^o) dt + Q(t)dW(t),~~ \psi(T) = \Phi_x(x^o(T)). \label{eq34}  \end{eqnarray}

Note the similarity in appearance with the Pontryiagin minimum principle. In fact we recover the  Pontryagin minimum principle for relaxed controls in \cite{ahmed2005,ahmed2006} by setting $\sigma =0.$  \\
For  controls based on  full-information which are ${\cal F}_t$ adapted, and under the condition that $\{{\cal F}_t, t  \in [0,T]\}$ is the natural filtration generated by the Brownian motion $\{W(t), t \in [0,T]\}$, augmented by all $P-$null sets in ${\cal F}$,   given by the inequality  (\ref{eq16}) (or equivalently (\ref{eq32})) is equivalent to the following point wise almost sure inequality (the derivation is similar to that of Corollary 5.3):

 \begin{eqnarray*}  H(t,x^o(t),\psi(t),Q(t),\mu )    \geq  H(t,x^o(t),\psi(t),Q(t),u_t^o), \nonumber \\
\:  \forall \mu \in {\cal M}_1(U),  \: a.e. \: t \in [0,T], \: P-a.s.    \label{eq35a}  \end{eqnarray*}
or equivalently,

 \begin{eqnarray*}  H(t,x^o(t),\psi(t),Q(t),u_t^o )   = \min_{ \mu \in {\cal M}_1(U)}   H(t,x^o(t),\psi(t),Q(t),\mu),  \\
\:    a.e. \: t \in [0,T], \: P-a.s.    \label{eq35ab}  \end{eqnarray*}
subject to the Hamiltonian system (\ref{eq33})-(\ref{eq34}).

 \vskip6pt  For the partial information case,   the point wise necessary conditions of optimality for controls are given in the next Corollary.

 \vskip6pt \noindent{\bf Corollary 5.3} Suppose the assumptions of Theorem 5.1 hold and consider controls which are ${\cal G}_t$ adapted.  Then the inequality  (\ref{eq16}) (or equivalently (\ref{eq32})) is equivalent to the following point wise almost sure inequality with respect to the $\sigma$-algebra ${\cal G}_t \subset {\cal F}_t:$

\begin{eqnarray}  {\cal E} \{  H(t,x^o(t),\psi(t),Q(t),\mu )|{\cal G}_t\} \geq   {\cal E}\{H(t,x^o(t),\psi(t),Q(t),u_t^o)|{\cal G}_t\} \label{eq35} \end{eqnarray}  for all $\mu \in {\cal M}_1(U)$, $a.e. t \in [0,T], P-a.s.$ subject to the Hamiltonian system (\ref{eq33})-(\ref{eq34}).

\vskip6pt\noindent{\bf Proof.} Since the admissible controls are  vaguely  ${\cal G}_t$ measurable, we can rewrite  the inequality (\ref{eq32}) in the following equivalent form,

\begin{eqnarray} &~& {\cal E} \int_0^T  {\cal E} \{  H(t,x^o(t),\psi(t),Q(t),u_t)|{\cal G}_t\}~  dt \nonumber \\ &~& ~~~~~~~~~~~~~ \geq   {\cal E} \int_0^T  {\cal E}\{H(t,x^o(t),\psi(t),Q(t),u_t^o)|{\cal G}_t \}~  dt. \label{eq36}  \end{eqnarray}

Let  $t \in (0,T),$  $\omega \in \Omega$ and $\varepsilon >0$ and  consider the sets $I_{\varepsilon} \equiv [t,t+\varepsilon] \subset I $ and  $\Omega_{\varepsilon} (\subset \Omega) \in {\cal G}_t$  containing $\omega$ such that $|I_{\varepsilon}| \rightarrow 0$  and $P(\Omega_{\varepsilon}) \rightarrow 0$ as $\varepsilon \rightarrow 0.$ For any subsigma algebra ${\cal G} \subset {\cal F}$, let $P_{{\cal G}}$ denote the restriction of the probability measure $P$ on to the $\sigma$-algebra ${\cal G}.$   For any (vaguely) ${\cal G}_t$-measurable  $\nu \in {\cal M}_1(U),$ construct the control  $$ u_t = \begin{cases}  \nu & \mbox{for}~~ (t,\omega) \in I_{\varepsilon} \times \Omega_{\varepsilon}  \\   u_t^o & \mbox{ otherwise}.\end{cases}.$$ Clearly, it follows from the above construction that $u \in {\cal U}_{ad}.$  Using this control in (\ref{eq35}) we obtain the following inequality
\begin{eqnarray} \int_{\Omega_{\varepsilon}\times I_{\varepsilon}} {\cal E} \{  H(t,x^o(t),\psi(t),Q(t),\nu)|{\cal G}_t\}~dt  \geq  \int_{\Omega_{\varepsilon}\times I_{\varepsilon}} {\cal E} \{  H(t,x^o(t),\psi(t),Q(t),u_t^o)|{\cal G}_t\}~ dt.\nonumber \\ &~& \label{eq37}  \end{eqnarray} Letting $|I_{\varepsilon}|$ denote the Lebesgue  measure of the set $I_{\varepsilon}$ and dividing the above expression  by the product measure $P(\Omega_{\varepsilon})|I_{\varepsilon}|$ and letting $\varepsilon \rightarrow 0$ we arrive at the following in equality,
\begin{eqnarray*}  {\cal E} \{  H(t,x^o(t),\psi(t),Q(t),\nu)|{\cal G}_t\}  \geq   {\cal E} \{  H(t,x^o(t),\psi(t),Q(t),u_t^o)|{\cal G}_t\}~ \end{eqnarray*} which holds  for almost all $t \in I$ and $P_{{\cal G}_t}$ almost all $\omega \in \Omega.$ Thus we have completed the proof. $\bullet$ \vskip6pt\noindent {\bf Remark 5.4}  Define $$ g_t(\xi) \equiv {\cal E}\{ H(t,x^o(t),\psi(t),Q(t),\xi)|{\cal G}_t \}, t \in I,  \xi \in U.$$ The reader can easily verify from the basic assumptions on the parameters $\{b,\sigma,\ell,\Phi\}$ that the random process $g$  is an element of $ L_1^a(I,C(U))$ and that it is adapted to the $\sigma$-algebra ${\cal G}_t.$ Clearly,   the necessary condition given by the  inequality (\ref{eq35})  can be written as follows
\begin{eqnarray*} \int_{U} g_t(\xi)  \mu(d\xi)  \geq  \int_{U} g_t(\xi) u_t^o(d\xi), \end{eqnarray*} and this  must hold  for all ${\cal M}_1(U)$-valued ${\cal G}_t$-adapted  (vaguely  measurable) random variables $\mu.$   Define  $$\Lambda_t(\mu) \equiv \int_{U} g_t(\xi) \mu(d\xi).$$  This is a ${\cal G}_t$- measurable continuous linear functional on ${\cal M}_1(U)$. Since the later space is vaguely compact, it attains its minimum on ${\cal M}_1(U)$ and from the above inequality it follows that $u_t^o$ is one such element.  Because the functional $\Lambda_t$ is not strictly convex there may be multiplicities of minima $M^o(t).$   It is easy to verify that the set $$ M^o(t) \equiv \{ \mu \in M_1(U): \mu~\hbox{is}~  {\cal G}_t-\hbox{measurable}~ \hbox{and}~  \Lambda_t(\mu) = \Lambda_t(u^o_t)\}$$ is convex and a  vaguely (weak star) closed subset of ${\cal M}_1(U)$ and hence vaguely  compact. Thus  $t \longrightarrow M^o(t)$ is a measurable multi function with convex compact   values in ${\cal M}_1(U).$ By our assumption $U$ is compact and hence ${\cal M}_1(U)$ is a compact Polish space and hence  a compact  Souslin space.  Thus it follows from the well known Yankov-Von Neumann-Auman selection theorem [\cite{Hu-Papageorgiou1997}, Theorem 2.14, p158]  that the multi function   $t \longrightarrow M^o(t)$ has a ${\cal G}_t$ measurable selection. Hence we   have a ${\cal G}_t$  measurable  optimal relaxed control.

\section{Extension to Jump Processes}
\label{jump}

The necessary conditions of optimality  given in the previous section can be easily extended to control problems involving stochastic differential equations driven both by Brownian motion and Le\'vy process or Poisson jump  process. Let $Z \equiv R^n \setminus\{0\}$ and  ${\cal B}(Z)$ the Borel algebra of subsets of the set $Z$. Let $p(dv\times dt)$ denote the Poisson counting measure on ${\cal B}(Z)\times \sigma(I).$   Physical interpretation of this measure is simple. For each  $\Gamma \in {\cal B}(Z)$ and any interval $\Delta \in \sigma(I)$,  $p(\Gamma \times \Delta)$ gives the number of  jumps over the interval $\Delta$ of sizes confined  in $\Gamma.$  This is a Poisson  random variable with mean  $ {\cal E} p(\Gamma \times \Delta) = \pi(\Gamma) \lambda(\Delta)$ where $\lambda$ is the Lebesgue measure on the real line and $\pi$ is the  Le\'vy measure on $Z.$  Here $\pi$  is a countably additive bounded positive measure.  The compensated Poisson random measure is given by $$ q(dv\times dt) = p(dv\times dt) - \pi(dv) dt.$$ There is no loss of generality considering the compensated Poission random measure  in modeling SDE. As usual, we assume that all the random processes considered in this paper  are based on the  filtered probability space $(\Omega,{\cal F}, {\cal F}_{t\geq 0},P)$ where $\{{\cal F}_t, t \geq 0\} $ is an increasing family of subsigma algebras of $\sigma$-algebra  ${\cal F}$ and that they are right continuous with left limits.
A controlled  stochastic differential equation driven both  by Brownian motion and the compensated jump process described above is given by the following stochastic differential equation \begin{eqnarray} dx(t) = b(t,x(t),u_t) dt + \sigma(t,x(t),u_t) dW(t) + \int_{Z} C(t,x(t),v,u_t) q(dv\times dt), t \in I \label{eq38} \end{eqnarray} for $ x(0) = x_0.$  Throughout the rest of the  paper it is assumed without any further notice  that $\{x_0, W,q\}$ are independent random elements. Again our controls are relaxed controls  which,  for the partial information case,   are weakly ${\cal G}_t$ adapted ${\cal M}_1(U)$ valued random processes denote by ${\cal U}_{ad}.$ The cost functional is given by \begin{eqnarray} J(u) \equiv {\cal E} \biggl\{ \int_I \ell(t,x(t),u_t) dt + \Phi(x(T))\biggr\}. \label{eq39}\end{eqnarray}  Objective is to find a control from the admissible set ${\cal U}_{ad}$  at which the functional (\ref{eq39}) attains its minimum.   The method of proof of the necessary conditions of optimality for this model is no different from the one given for  the continuous case.  Hence we present the results without repeating the detailed proof.

For the problem involving jump process, we  introduce the following  Hilbert space of discontinuous square integrable  semi martingales denoted by ${\cal DSM}_0^2$ and this is given by
\begin{eqnarray} \hskip-20pt  &~& {\cal DSM}_0^2 \equiv \biggl \{ m: m(t) = \int_0^t v(s) ds + \int_0^t Q(s) dW(s) + \int_0^t \int_{Z} \varphi(v,t) q(dv\times dt) \nonumber \\ &~& ~~~~~~~~~~~~~~~v \in L_2^a(I,R^n), Q \in L_2^a(I,{\cal L}(R^m,R^n)), \varphi \in L_2^a(I,L_2^n(Z,\pi))  \biggr\} \label{eq40} \end{eqnarray} where $L_2^n(Z,\pi)$ denotes the Hilbert space of $R^n$-valued  functions defined on $Z$ which are square integrable with respect to the Le\'vy measure $\pi.$    In this case the  norm topology is  given  by \begin{eqnarray} \parallel m\parallel_{{\cal DSM}_0^2} = \biggl( {\cal E}\int_I |v(t)|_{R^n}^2 dt + {\cal E} \int_{I} tr(Q^*(t)Q(t)) dt + {\cal E}  \int_I \int_{Z} |\varphi(v,t)|_{R^n}^2 \pi(dv) dt \biggr)^{1/2}. \label{eq41} \end{eqnarray}
Now we are prepared to present  the necessary conditions of optimality. Before we  do so we need the following assumptions for $C.$ \vskip6pt  The function   $C: I\times R^n\times R^n \times U \longrightarrow  R^n $ is measurable in $t$ on $I$ and continuous in the rest of the arguments  satisfying, uniformly with respect to $\xi \in U,$ the following assumptions   \par \begin{eqnarray*} (A6): ~~~~~~~~~~~~~~~~~~~~~~~~ \biggl( \int_{Z} |C(t,x,v,\xi)|_{R^n}^2\pi(dv)\biggr)^{1/2} \leq K(t) (1 + |x|_{R^n}) \end{eqnarray*}
\par \begin{eqnarray*} (A7):~~~~~~~ \biggl( \int_{Z} |C(t,x,v,\xi)-C(t,y,v,\xi)|_{R^n}^2\pi(dv)\biggr)^{1/2} \leq K(t) (|x-y|_{R^n}). \end{eqnarray*}

\vskip6pt\noindent {\bf Theorem 6.1} Consider the system (\ref{eq38}) with the cost functional (\ref{eq39}) and the admissible controls ${\cal U}_{ad}.$ Suppose $\{b,\sigma,C\}$ satisfy the assumptions (A1)-(A7) and that their first derivatives with respect to the state variable $x \in R^n$ are  uniformly bounded.   An element $ u^o \in {\cal U}_{ad},$ with the corresponding solution $x^o \in B_{\infty}^a(I, L_2(\Omega,R^n))$ to be optimal, it is necessary that there exists a semi  martingale  $m^o \in {\cal DSM}_0^2$ with the intensity process $(\psi,Q,\varphi) \in  L_2^a(I,R^n)\times L_2^a(I,{\cal L}(R^m,R^n))\times L_2^a(I,L_2^n(Z,\pi))$ such that the following inequality and the stochastic differential equations  hold:

\begin{eqnarray} (1): \hskip-20pt &~& {\cal E} \int_0^T \biggl\{ (b(t,x^o(t),u_t-u_t^o),\psi(t)) + tr(Q^*(t)\sigma(t,x^o(t),u_t-u_t^o)) \nonumber \\ &~& ~~+ \int_{Z} (C(t,x^o(t),v,u_t-u_t^o),\varphi(t,v))\pi(dv) +  \ell(t,x^o(t),u_t-u_t^o)\biggr\} dt  \geq 0 \label{eq42}\end{eqnarray} for all $u \in {\cal U}_{ad}.$
\begin{eqnarray} (2): \hskip-20pt  &~&    dx^o(t) = b(t,x^o(t),u_t^o) dt + \sigma(t,x^o(t),u_t^o)dW(t) + \int_{Z} C(t,x^o(t),v,u_t^o) q(dv\times dt)  \nonumber \\ &~& x^o(0) = x_0 \label{eq43}  \end{eqnarray}

 \begin{eqnarray} (3):  \hskip-20pt  &~&   -d \psi(t)  = b_x^*(t,x^o(t),u_t^o)\psi(t) dt +  V_{Q}(t) dt  - Q(t) dW(t)  \nonumber \\
~~~~~~~~ &+& \int_{Z} C_x^*(t,x^o(t),v,u_t^o)\varphi(t,v) \pi(dv) dt  - \int_{Z} \varphi(t,v) q(dv\times dt)  +\ell_x(t,x^o(t),u_t^o) dt  \nonumber \\ &~& \psi(T) = \Phi_x(x^o(T)) \label{eq44}  \end{eqnarray} where  $V_{Q} \in L_2^a(I,R^n)$ is   given by  $(V_{Q}(t),\zeta) = tr (Q^*(t)\sigma_x(t,x^o(t),u_t^o; \zeta)), t \in I.$
\vskip6pt \noindent  {\bf Remark 6.2}  Define the Hamiltonian $$H : I \times R^n\times R^n\times {\cal L}(R^m,R^n)\times L_2^n(Z,\pi) \times {\cal M}_1(U) \longrightarrow R $$ by  the following expression \begin{eqnarray}&~&  H(t,x,\psi,Q,\varphi,\mu) \equiv (b(t,x,\mu),\psi) + tr (Q^*\sigma(t,x,\mu)) \nonumber \\ &~& ~~~~~~~~~~~~~~~~~~~~~~~+  \int_{Z} (C(t,x,v,\mu),\varphi(t,v))_{R^n}\pi(dv) + \ell(t,x,\mu), \label{eq45}\end{eqnarray} where $\varphi \in L_2^n(Z,\pi).$ We write $t\longrightarrow \varphi(t)$ for the $L_2^n(Z,\pi)$ valued function.  In terms of this Hamiltonian, the necessary conditions of Theorem 6.1 can be written in the following canonical form:
\begin{eqnarray}&~&  {\cal E} \int_I H(t,x^o(t),\psi(t),Q(t),\varphi(t),u_t) dt  \geq {\cal E} \int_I H(t,x^o(t),\psi(t),Q(t),\varphi(t),u_t^o) dt\nonumber \\ &~& ~\hbox{ for all}~ u \in {\cal U}_{ad}, \label{eq46} \\ &~& dx^o(t) = H_{\psi} dt + \sigma(t,x^o(t),u_t^o) dW(t) \nonumber \\
~~~&~& + \int_{Z} C(t,x^o(t),v,u_t) q(dv\times dt),~~ x^o(0) = x_0 \label{eq47} \\ &~& d \psi(t) =  -H_x dt + Q(t) dW(t) + \int_{Z} \varphi(t,v) q(dv\times dt),~~ \psi(T) = \Phi_x(x^o(T)). \label{eq48} \end{eqnarray}

Similarly as before, one can also obtain point wise almost sure variational inequalities.

\section{Necessary conditions with Regular Controls}
\label{regular}

In the development of the necessary conditions of optimality given in the preceding two sections we have tacitly  used the existence Theorem 3.2  which asserts the existence of optimal controls from the class of relaxed controls ${\cal U}_{ad}.$  Let $ L_{\infty}^a(I\times \Omega, U)$ denote the class of ${\cal G}_t$ adapted random processes defined on the interval $I$ and taking values from the closed bounded set $U\subset R^d.$  This is the class of regular controls and we denote this by ${\cal U}^{r}.$ It is clear that this embeds continuously  into the class of relaxed controls through the map $u \ni {\cal U}^r \longrightarrow \delta_{u(t,\omega)}\in {\cal U}_{ad}.$   Clearly, for every $\vartheta \in L_1^a(I\times \Omega, C(U))$ \begin{eqnarray}  {\cal E} \int_{I\times U}\vartheta(t,\omega,\xi) \delta_{u(t,\omega)}(d\xi) dt  = {\cal E}\int_{I} \vartheta (t,\omega,u(t,\omega))dt. \label{eq49} \end{eqnarray}

\vskip6pt\noindent {\bf Theorem 7.1}   Consider the class of regular controls ${\cal U}^r$  with $U$ assumed to be closed bounded and convex. Suppose Theorem 3.2 holds for regular controls  in the sense that an optimal control exists from the class ${\cal U}^r.$  Then  all the necessary conditions involving relaxed controls (Theorem 5.1, Theorem 6.1)  reduce to the classical minimum principle for stochastic systems. \vskip6pt \noindent {\bf Proof.} The proof is direct.  In fact it follows from straightforward application of the embedding mentioned above and the definition  (\ref{eq49}). Considering  the necessary conditions of optimality given by Theorem 5.1, and using the  embedding mentioned above it is easy to derive the following necessary conditions of optimality \begin{eqnarray} (1): \hskip-20pt &~& {\cal E} \int_0^T \bigl\{ (b(t,x^o(t),u_t),\psi(t)) + tr(Q^*(t)\sigma(t,x^o(t),u_t)) + \ell(t,x^o(t),u_t)\bigr\} dt \nonumber \\ &~& ~~~~\geq {\cal E} \int_0^T \bigl\{ (b(t,x^o(t),u_t^o),\psi(t)) + tr(Q^*(t)\sigma(t,x^o(t),u_t^o)) + \ell(t,x^o(t),u_t^o)\bigr\} dt,  \nonumber  \\ &~& ~~\hbox{ for all }~~ u \in {\cal U}^r. \label{eq50}\end{eqnarray}

\begin{eqnarray} (2): \hskip-20pt  &~&    dx^o(t) = b(t,x^o(t),u_t^o) dt + \sigma(t,x^o(t),u_t^o))dW(t) \nonumber \\ &~& x^o(0) = x_0 \label{eq51}  \end{eqnarray}

 \begin{eqnarray} (3):  \hskip-20pt &~&    -d \psi(t)  = b_x^*(t,x^o(t),u_t^o)\psi(t)  dt + V_{Q}(t) dt +\ell_x(t,x^o(t),u_t^o) dt - Q(t) dW(t)  \nonumber \\ &~& \psi(T) = \Phi_x(x^o(T)) \label{eq52}  \end{eqnarray} where  $V_{Q} \in L_2^a(I,R^n)$ is   given by  $(V_{Q}(t),\zeta) = tr (Q^*(t)\sigma_x(t,x^o(t),u^o(t); \zeta)), t \in I.$   \vskip6pt \noindent {\bf Remark 7.2} Using precisely similar  arguments for the SDE with jumps,    one can obtain  the minimum principle for regular controls from those of relaxed controls  given by Theorem 6.1. \section{Realizability  of Relaxed Controls by Regular Controls}

  We proved existence of optimal relaxed controls in Theorem 3.2 without requiring convexity of the control domain $U$. In any application it is much easier to construct regular controls. So one may be interested to find a regular control corresponding to which the performance of the system is close to that  realized by optimal  relaxed control. In this regard we have the following result.

  \vskip6pt\noindent {\bf Theorem 8.1} Consider the regular controls ${\cal U}^r$ with $U$ closed bounded but not necessarily convex as in Theorem 7.1.  Suppose the basic assumptions of Lemma 3.1 and Theorem 3.2 hold  and  consider the control problem as stated in  Theorem 3.2. Further, suppose that  $ x \longrightarrow \Phi(x)$ is continuous.  Let  $u^o \in {\cal U}_{ad}$ be  the optimal relaxed control. Then,  for every  $\varepsilon >0$ there exists a regular control $u_r \in {\cal U}^r$ such that $$ J(u_r) \leq \varepsilon + J(u^o).$$

 \vskip6pt\noindent {\bf Proof} Since ${\cal U}_{ad} \equiv L_{\infty}^a(I\times \Omega, {\cal M}_1(U)) \subset  L_{\infty}^a(I\times \Omega, {\cal M}(U)) $ is  compact in the vague topology (that is weak star topology) and convex (because ${\cal M}_1(U)$ is convex, it follows from the  well  known Krein-Millman theorem that $$ {\cal U}_{ad} = cl^{v}conv (\hbox{ext}( {\cal U}_{ad})),$$ that is, ${\cal U}_{ad}$  is the weak star closed convex hull of its extreme points. Considering the embedding ${\cal U}^r \hookrightarrow  {\cal U}_{ad}$ as  mentioned above,  it is easy to verify that the extreme points of ${\cal U}_{ad}$ are precisely  the set of regular controls ${\cal U}^r$ through the map $u \ni {\cal U }^r \longrightarrow \delta_{u} \in {\cal U}_{ad}.$ Thus, if $u^o \in {\cal U}_{ad}$ is the optimal (relaxed) control there exists a sequence $\{u^n\}$ of the form  $$u^n \equiv \sum_{i=1}^n \alpha_i^n u_i, u_i \in {\cal U}^r, \alpha_i^n \geq 0, \sum_{i=1}^n \alpha_i^n =1, n\in N $$  such that $u^n \buildrel v \over\longrightarrow u^o.$ Let $\{x^n, x^o\} \subset B_{\infty}^a(I,L_2(\Omega,R^n))$ denote the solutions of the system equation (\ref{eq1}) corresponding to the  controls $\{u^n,u^o\}$ respectively.  Then it follows from Lemma 3.1 that, along a subsequence if necessary, $x^n \buildrel s \over\longrightarrow x^o$ in $B_{\infty}^a(I,L_2(\Omega,R^n)).$ Consequently,  it follows from continuity of $\ell$ and $\Phi$  in the state variable $x$ and the assumptions (a1)-(a3) and Lebesgue dominated convergence theorem that $\lim_{n\rightarrow \infty} J(u^n) = J(u^o).$ Note that for every $n \in N$, $u^n \in {\cal U}^r$, and so, for every $\varepsilon >0$, there exists an $ n_{\varepsilon}\in N$ such that $ |J(u^n)- J(u^o)| < \varepsilon $ for all $n \geq n_{\varepsilon}.$ Taking $u_r = u^{n_{\varepsilon}}$ we have $ J(u_r) \leq \varepsilon + J(u^o).$  This completes the proof. $\bullet$
\vskip6pt\noindent {\bf Remark 8.2} In view of the above result it is evident that an $\epsilon$-optimal control can be found from the class of regular controls (measurable functions with values in $U$) though the limit of such controls may be  a relaxed control.  More specifically if $U \subset R^d$ consists of a finite set of points, it is clearly non-convex,  and optimal control may not exist from the class of regular controls ${\cal U}^r$ based on the set $U.$ However, optimal relaxed controls do exist. In this case the sequence of regular controls approximating the optimal relaxed control may oscillate violently between the finite set of points of $U$ with increasing frequency (converging to infinity). This  is known as  chattering.

\bibliographystyle{IEEEtran}
\bibliography{bibdata}

\begin{thebibliography}{10}
\providecommand{\url}[1]{#1}
\csname url@samestyle\endcsname
\providecommand{\newblock}{\relax}
\providecommand{\bibinfo}[2]{#2}
\providecommand{\BIBentrySTDinterwordspacing}{\spaceskip=0pt\relax}
\providecommand{\BIBentryALTinterwordstretchfactor}{4}
\providecommand{\BIBentryALTinterwordspacing}{\spaceskip=\fontdimen2\font plus
\BIBentryALTinterwordstretchfactor\fontdimen3\font minus
  \fontdimen4\font\relax}
\providecommand{\BIBforeignlanguage}[2]{{%
\expandafter\ifx\csname l@#1\endcsname\relax
\typeout{** WARNING: IEEEtran.bst: No hyphenation pattern has been}%
\typeout{** loaded for the language `#1'. Using the pattern for}%
\typeout{** the default language instead.}%
\else
\language=\csname l@#1\endcsname
\fi
#2}}
\providecommand{\BIBdecl}{\relax}
\BIBdecl

\bibitem{cesari1983}
L.~Cesari, \emph{Optimization Theory and Applications}.\hskip 1em plus 0.5em
  minus 0.4em\relax Springer-Verlag, 1983.

\bibitem{ahmed2006}
N.~U. Ahmed, \emph{Dynamic Systems and Control with Applications}.\hskip 1em
  plus 0.5em minus 0.4em\relax World Scientific, New Jersey London, Singapore
  Beijing Shanghai, Hong Kong, Taipei, Chenna, 2006.

\bibitem{ahmed1981}
N.~U. Ahmed and K.~L. Teo, \emph{Optimal Control of Distributed Parameter
  Systems}.\hskip 1em plus 0.5em minus 0.4em\relax Elsevier North Holland, New
  York, Oxford, 1981.

\bibitem{Fattorini1999}
H.~O. Fattorini, \emph{Infinite Dimensional Optimization and Control
  Theory}.\hskip 1em plus 0.5em minus 0.4em\relax Encyclopedia of Mathematics
  and Its Applications, 62, Cambridge University Press., 1999.

\bibitem{ahmed2005}
N.~U. Ahmed, ``Optimal relaxed controls for systems governed by impulsive
  differential inclusions,'' \emph{Nonlinear Functional Analysis \&
  Applications}, vol.~10, no.~3, pp. 427--460, 2005.

\bibitem{ahmed-charalambous2007}
N.~U. Ahmed and C.~D. Charalambous, ``Minimax games for stochastic systems
  subject to relative entropy uncertainty: Applications to {SDE}'s on {H}ilbert
  spaces,'' \emph{Journal of Mathematics of Control, Signals and Systems},
  vol.~19, pp. 65--91, 2001.

\bibitem{kushner1972}
H.~J. Kushner, ``Necessary conditions for continuous parameter stochastic
  optimization problems,'' \emph{SIAM Journal on Control and Optimization,},
  1972.

\bibitem{haussmann1986}
U.~Haussmann, \emph{A Stochastic Maximum Principle for Optimal Control of
  Diffusions}, ser. Longman Sci. \& Tech., Harlow, UK.\hskip 1em plus 0.5em
  minus 0.4em\relax Pitman Research Notes in Mathematics, 1986, vol. 151.

\bibitem{bismut1978}
J.~Bismut, ``An introductory approach to duality in optimal stochastic
  contro,'' \emph{SIAM Review}, vol.~30, pp. 62--78, 1978.

\bibitem{bensoussan1981}
A.~Bensoussan, \emph{Lecture on Stochastic Control}, ser. Lecture Notes in
  Mathematics.\hskip 1em plus 0.5em minus 0.4em\relax Springer-Verlag, Berlin,
  1982.

\bibitem{bensoussan1983}
------, ``Stochastic maximum principle for distributed parameter systems,''
  \emph{Journal of the Franklin Institute}, vol. 315, pp. 387--406, 1999.

\bibitem{peng1990}
S.~Peng, ``A general stochastic maximum principle for optimal control
  problems,'' \emph{SIAM Journal on Control and Optimization}, vol.~28, no.~4,
  pp. 966--979, 1990.

\bibitem{cadenillas-karatzas1995}
A.~Cadenillas and I.~Karatzas, ``The stochastic maximum principle for linear
  convex optimal control with random coefficients,'' \emph{SIAM Journal on
  Control and Optimization}, 1995.

\bibitem{elliott-kohlmann1994}
R.~Elliott and M.~Kohlmann, ``The second order minimum principle and adjoint
  process,'' \emph{Stochastics \& Stochastic Reports,}, vol.~46, pp. 25--39,
  1994.

\bibitem{zhou1991}
X.~Y. Zhou, ``A unified treatment of maximum principle and dynamic programming
  in stochastic controls,'' \emph{Stochastics \& Stochastic Reports}, vol.~36,
  pp. 137--161, 1991.

\bibitem{yong-zhou1999}
J.~Yong and X.~Zhou, \emph{Stochastic Controls, Hamiltonian Systems and HJB
  Equations}.\hskip 1em plus 0.5em minus 0.4em\relax Springer-Verlag, 1999.

\bibitem{mezerdi-bahlali2002}
B.~Mezerdi and S.~Bahlali, ``Necessary conditions for optimality in relaxed
  stochastic control problems,'' \emph{Stochastics \& Stochastic Reports},
  vol.~73, pp. 201--218, 2002.

\bibitem{bahlali-mezerdi-djehiche2006}
B.~M. S.~Bahlali and B.~Djehiche, ``Approximation and optimality necessary
  conditions in relaxed stochastic control problems,'' \emph{Journal of Applied
  Mathematics and Stochastic Analysis}, p.~23, 2006, iD 72762.

\bibitem{bahlali-djehiche-mezerdi2007a}
B.~D. S.~Bahlali and B.~Mezerdi, ``The relaxed stochastic maximum principle in
  singular optimal control of diffusions,'' \emph{SIAM Journal on Control and
  Optimization}, vol.~46, pp. 427--444, 2007.

\bibitem{bahlali-djehiche-mezerdi2007b}
K.~D. S.~Bahlali and B.~Mezerdi, ``On the stochastic maximum principle in
  optimal control of degenerate diffusions with lipschitz coefficients,''
  \emph{Applied Mathematics and Optimization}, vol.~56, pp. 364--378, 2007.

\bibitem{bahlali2008}
S.~Bahlali, ``Necessary and sufficient optimality coditions for relaxed and
  strict control problems,'' \emph{SIAM Journal on Control and Optimization},
  2008.

\bibitem{buckdahn-djehiche-li2011}
B.~D. R.~Buckdahn and J.~Li, ``A general stochastic maximum principle for sde's
  of mean-field type,'' \emph{Applied Mathematics and Optimization}, vol.~64,
  no. 197-216, 2011.

\bibitem{zhang-shi2011}
L.~Zhang and Y.~Shi, ``Maximum principle for forward-backward doubly stochastic
  control systems and applications,'' \emph{ESAIM: Control,Optimization and
  Calculus of Variations}, vol. COCV 17, pp. 1174--1197, 2011, dOI:
  10,1051/cocv/2010042.

\bibitem{zhang-elliott-siu2012}
R.~J.~E. X.~Zhang and T.~K. Siu, ``A stochastic maximum principle for a markov
  regime-switching jump-diffusion model and its applications to finance,''
  \emph{SIAM Journal on Control and Optimization}, vol.~50, no.~2, pp.
  964--990, 2012.

\bibitem{charalambous-hibey1996}
C.~D. Charalambous and J.~L. Hibey, ``Minimum principle for partially
  observable nonlinear risk-sensitive control problems using measure-valued
  decompositions,'' \emph{Stochastics \& Stochastic Reports}, 1996.

\bibitem{Hu-Papageorgiou1997}
S.~Hu and N.~S. Papageorgiou, \emph{Handbook of Multivalued Analysis}.\hskip
  1em plus 0.5em minus 0.4em\relax Kluwer Academic Publishers, Dordrecht,
  Boston, London., 1997.

\end{thebibliography}

\end{document}